\documentclass[11pt]{article}

\newtheorem{theorem}{Theorem}
\newtheorem{definition}{Definition}

\usepackage{amssymb}
\usepackage{amsmath}
\usepackage{float}
\def\ds{\displaystyle}
\def\dt{k}
\def\dx{h}
\def\eps{\varepsilon}
\def\rmr{{\mathrm{r}}}
\def\rms{{\mathrm{s}}}
\def\bfB{{\mathbf{B}}}
\def\bfD{{\mathbf{D}}}
\def\bfE{{\mathbf{E}}}
\def\bfJ{{\mathbf{J}}}
\def\bfP{{\mathbf{P}}}
\def\bfj{{\mathbf{j}}}
\def\bfx{{\mathbf{x}}}
\def\bfxi{{\boldsymbol{\xi}}}
\def\bbC{{\mathbb{C}}}
\def\bbN{{\mathbb{N}}}
\def\bbR{{\mathbb{R}}}
\def\curl{{\operatorname{\ curl\ }}}
\def\transp{{}^t}
\def\uv{\,\textrm{m\,s}^{-1}}
\def\ut{\,\textrm{s}}
\def\um{\,\textrm{m}}
\def\uf{\,\textrm{rad\,s}^{-1}}
\begin{document}

\title{Stability of FD--TD schemes for Maxwell--Debye and Maxwell--Lorentz 
equations.}

\author{Brigitte Bid\'egaray-Fesquet\thanks{B. Bid\'egaray-Fesquet is with the
LMC-IMAG, CNRS UMR 5523, B.P. 53, 38041 Grenoble Cedex 9, France. 
E-mail: brigitte.bidegaray@imag.fr .}}

\maketitle

\begin{abstract}
The stability of five finite difference--time domain (FD--TD) sche\-mes 
coupling Maxwell equations to Debye or Lorentz models have been analyzed in 
\cite{petropoulos1}, where numerical evidence for specific media have been 
used. We use von Neumann analysis to give necessary and sufficient stability 
conditions for these schemes for any medium, in accordance with the partial 
results of \cite{petropoulos1}.
\end{abstract}


\section{Introduction}

To describe the propagation of an electromagnetic wave through a 
dispersive medium some extensions to Maxwell equations are used. They involve 
time differential equations which accounts for the constitutive laws of the 
material that link the displacement $\bfD$ to the electric field $\bfE$ or 
equivalently the polarization $\bfP$ to $\bfE$. We focus on two of these 
models (Debye and Lorentz models) which are addressed in \cite{petropoulos1} 
in view of specific applications to the interaction of an electromagnetic 
wave with a human body. In contrast we treat any medium which is described by 
these models. We only consider the stability analysis of numerical schemes 
whereas \cite{petropoulos1} also treated phase error issues.

\subsection{Maxwell--Debye and Maxwell--Lorentz models}

In our context (no magnetization) the Maxwell equations read
\begin{equation}
\label{Maxwell}
\begin{array}{lrcl}
\textrm{(Faraday)} & 
\ds\partial_t \bfB(t,\bfx) & = & \ds- \curl \bfE(t,\bfx), \\
\textrm{(Amp\`ere)} &
\ds\partial_t \bfD(t,\bfx) & = & \ds\frac1{\mu_0} \curl \bfB(t,\bfx),
\end{array}
\end{equation}
where $\bfx\in\bbR^N$ together with a linear constitutive law
\begin{equation}
\label{constitutive}
\bfD(t,\bfx) = \eps_0 \eps_\infty\bfE(t,\bfx) 
+ \eps_0 \int_{-\infty}^t \bfE(t-\tau,\bfx) \chi(\tau)d\tau,
\end{equation}
where $\eps_\infty$ is the relative infinite frequency permittivity and $\chi$
is the linear susceptibility. The discretization of the integral expression 
\eqref{constitutive} leads to recursive schemes (see e.g. 
\cite{luebbers-hunsberger-kunz-standler-schneider},
\cite{young-kittichartphayak-kwok-sullivan}). However, differentiating Eq.
\eqref{constitutive} leads to a time differential equation for $\bfD$ which 
depends on the specific form of $\chi$. For a Debye medium 
\begin{equation}
\label{DebyeD}
t_\rmr \partial_t \bfD + \bfD = t_\rmr \eps_0 \eps_\infty \partial_t \bfE 
+ \eps_0\eps_\rms \bfE,
\end{equation}
where $t_\rmr>0$ is the relaxation time and $\eps_\rms\geq\eps_\infty$ is the 
relative static permittivity. Defining the polarization by 
$\bfP(t,\bfx)= \bfD(t,\bfx)-\eps_0 \eps_\infty\bfE(t,\bfx)$, an equivalent 
form is  
\begin{equation}
\label{DebyeP}
t_\rmr \partial_t \bfP + \bfP = \eps_0(\eps_\rms-\eps_\infty) \bfE.
\end{equation}
For a Lorentz medium with one resonant frequency $\omega_1$, we likewise have
\begin{equation}
\label{LorentzD}
\partial_t^2 \bfD + \nu \partial_t \bfD + \omega_1^2\bfD 
= \eps_0\eps_\infty\partial_t^2 \bfE + \eps_0\eps_\infty\nu \partial_t \bfE 
+ \eps_0\eps_\rms\omega_1^2\bfE,
\end{equation}
where $\nu\geq0$ is a damping coefficient, and
\begin{equation}
\label{LorentzP}
\partial_t^2 \bfP + \nu \partial_t \bfP + \omega_1^2\bfP 
= \eps_0(\eps_\rms-\eps_\infty)\omega_1^2\bfE.
\end{equation}
If we denote by $\bfJ$ the time derivative of $\bfP$, system \eqref{Maxwell} 
can be cast as
\begin{equation}
\label{MaxwellP}
\begin{array}{rcl}
\ds\partial_t \bfB(t,\bfx) & = & \ds- \curl \bfE(t,\bfx), \\
\ds \eps_0\eps_\infty \partial_t \bfE(t,\bfx) 
& = & \ds \frac1{\mu_0} \curl \bfB(t,\bfx) - \bfJ(t,\bfx).
\end{array}
\end{equation}

\subsection{Numerical schemes}

A classical and very efficient way to compute the Maxwell equations is the Yee
scheme \cite{yee}. We restrict our study to existing Yee based schemes. Other
methods may be found in the literature in the context of Maxwell-Debye and 
Maxwell-Lorentz equations: see e.g. \cite{feise-schneider-bevelacqua} for 
pseudo-spectral schemes or \cite{stoykov-kuiken-lowery-taflove} for finite 
element--time domain (FE--TD) schemes.

The Yee scheme consists in discretizing $\bfE$ and $\bfB$ on staggered grids 
in space and time. This allows to use only centered discrete differential
operators. We denote by $\dx$ the space step (supposed here to be the same in 
all directions in the case of multi-dimensional equations) and by $\dt$ the 
time step.  In space dimension 1, we only consider the dependence in the space 
variable $z$ and classically two polarizations for the field may be decoupled. 
For example, the transverse electric polarization only involves $E\equiv E_x$ 
and $B\equiv B_y$. The discretized variables are $E_j^n \simeq E(n\dt,j\dx)$ 
(and similar notations for $D\equiv D_x$) and 
$B_{j+\frac12}^{n+\frac12} \simeq B((n+\frac12)\dt,(j+\frac12)\dx)$, and the 
Yee scheme for system \eqref{Maxwell} reads
\begin{equation}
\label{Max}
\begin{array}{rcl}
\ds \frac1\dt (B_{j+\frac12}^{n+\frac12}-B_{j+\frac12}^{n-\frac12})
& = & \ds - \frac1\dx (E_{j+1}^n-E_j^n), \\
\ds \frac1\dt (D_j^{n+1}-D_j^n) & = & \ds - \frac1{\mu_0 \dx}  
(B_{j+\frac12}^{n+\frac12}-B_{j-\frac12}^{n+\frac12}). 
\end{array}
\end{equation}
Similarly the Yee scheme for system \eqref{MaxwellP} reads
\begin{equation}
\label{MaxP}
\begin{array}{rcl}
\ds \frac1\dt (B_{j+\frac12}^{n+\frac12}-B_{j+\frac12}^{n-\frac12})
& = & \ds - \frac1\dx (E_{j+1}^n-E_j^n), \\
\ds \frac{\eps_0\eps_\infty}\dt (E_j^{n+1}-E_j^n) & = & 
\ds - \frac1{\mu_0 \dx} 
(B_{j+\frac12}^{n+\frac12}-B_{j-\frac12}^{n+\frac12})
- J_j^{n+\frac12}. 
\end{array}
\end{equation}

Usual Maxwell equations consist in taking $J_j^{n+\frac12}\equiv0$ in 
Eq.~\eqref{MaxP} or equivalently $D_j^n=\eps_0\eps_\infty E_j^n$ in 
Eq.~\eqref{Max} and leads to a stable second order scheme under a 
Courant--Friedrichs--Lewy (CFL) stability condition. Namely, if 
$c_\infty = 1/\sqrt{\eps_0\eps_\infty\mu_0}$ denotes the infinite frequency 
light speed, the CFL condition reads $c_\infty \dt\leq\dx$ if the space 
dimension is $N=1$ and $c_\infty \dt\leq\dx/\sqrt2$ for $N=2$ or 3.

In contrast to the recursive schemes, we are interested in direct integration 
schemes which are based on the finite difference--time domain (FD--TD) 
discretization of Eqs~\eqref{DebyeD} to \eqref{LorentzP} (see 
\cite{kashiwa-yoshida-fukai}, \cite{joseph-hagness-taflove}, \cite{young}).

\subsection{Outline}

The von Neumann stability analysis is recalled in Sect.~\ref{sec-vonNeumann}.
We also describe the sketch of our proofs which is common for all the schemes.
In Section~\ref{sec-Debye} two one dimensional direct integration schemes for 
Debye media are presented and analyzed, pointing carefully out the physical 
properties needed to ensure stability and the specific cases which have to be 
handled separately. Numerical applications to physical media are also given. 
The same point of view is carried out for Lorentz media in 
Section~\ref{sec-Lorentz}. Two-dimensional results are given in 
Section~\ref{sec-bidim}.

\section{Principles of the von Neumann analysis}
\label{sec-vonNeumann}

The von Neumann analysis allows to localize roots of certain 
classes of polynomials, which proves to be crucial here. We recall the main 
principles of this technique. Details and proofs of theorems may be found in 
\cite{strikwerda}.

\subsection{Schur and von Neumann polynomials}

We define two families of polynomials: Schur polynomials and simple von 
Neumann polynomials. 
\begin{definition}
A polynomial is a Schur polynomial if all its roots, $r$, satisfy $|r|<1$.
\end{definition}
\begin{definition}
A polynomial is a simple von Neumann polynomial if all its roots, $r$, lie on 
the unit disk ($|r|\leq1$) and its roots on the unit circle are simple roots.
\end{definition}
If a polynomial is of high degree or has sophisticated coefficients, it may be
difficult to locate its roots. However, there is a way to split this difficult
problem into many simpler ones. For this aim, we construct a sequence of 
polynomials of decreasing degree. Let $\phi$ be written as
\[
\phi(z) = c_0 + c_1 z + \dots + c_p z^p,
\]
where $c_0$, $c_1$ \dots, $c_p\in\bbC$ and $c_p\neq0$. We define its conjugate
polynomial $\phi^*$ by
\[
\phi^*(z) = c_p^* + c_{p-1}^* z + \dots + c_0^* z^p.
\]
Given a polynomial $\phi_0$, we may define a sequence of polynomials
\[
\phi_{m+1}(z) = \frac{\phi_m^*(0)\phi_m(z)-\phi_m(0)\phi_m^*(z)}{z}.
\]
It is clear that $\textrm{deg} \phi_{m+1} < \textrm{deg} \phi_m$, if 
$\phi_m\not\equiv 0$. Besides, we have the two following theorems.

\begin{theorem}
\label{Th_Schur}
A polynomial $\phi_m$ is a Schur polynomial of exact degree $d$ if and only 
if $\phi_{m+1}$ is a Schur polynomial of exact degree $d-1$ and 
$|\phi_m(0)|\leq|\phi_m^*(0)|$. 
\end{theorem}

\begin{theorem}
\label{Th_vonNeumann}
A polynomial $\phi_m$ is a simple von Neumann polynomial if and only if \\
\hspace*{1cm} $\bullet$ $\phi_{m+1}$ is a simple von Neumann polynomial and 
$|\phi_m(0)|\leq|\phi_m^*(0)|$, \\
or \\
\hspace*{1cm} $\bullet$ $\phi_{m+1}$ is identically zero and $\phi'_m$ is a 
Schur polynomial.
\end{theorem}

The main ingredient in the proof of both theorems is the Rouch\'e theorem 
(see \cite{strikwerda}). To analyze $\phi_0$, at each step $m$, conditions 
should be checked (leading coefficient is non-zero, 
$|\phi_m(0)|\leq|\phi_m^*(0)|$, \dots) until a definitive negative answer 
arises or the degree is 1. 

\subsection{Stability analysis}

The models we deal with are linear models. They may therefore be analyzed in
the frequency domain. Thus we assume that the scheme handles a variable 
$U^n_\bfj$ with spatial dependence 
\[
U^n_\bfj = U^n \exp(i\bfxi\cdot\bfj),
\]
where $\bfxi$ and $\bfj\in\bbR^N$, $N=1,2,3$. The amplification matrix $G$ is
the matrix such that $U^{n+1}=GU^n$. We assume that $G$ does not depend on 
time or on $\dx$ and $\dt$ separately  but only on the ratio $\dx/\dt$. Let 
$\phi_0$ be the characteristic polynomial of $G$, then we have a sufficient 
stability condition.
\begin{theorem}
\label{Th_CSstab}
A sufficient stability condition is that $\phi_0$ be a simple von Neumann 
polynomial.
\end{theorem}
This condition is not necessary. A scheme is stable if and only if the sequence
$(U^n)_{n\in\bbN}$ is bounded. Since we assume that $G$ does not depend on 
time, $U^n=G^nU^0$ and stability is also the boundedness of $(G^n)_{n\in\bbN}$.
If the eigenvalues of $G$, i.e. the roots $r$ of $\phi_0$, lie inside the unit
circle ($|r|<1$), then $\lim_{n\to\infty}G^n=0$ and the sequence is bounded. 
If any root lies outside the unit circle then $G^n$ grows exponentially and 
the scheme is unstable. The intermediate case when some roots may be on the 
unit circle (and the others inside) may lead to different situations. The good 
case is for example given when $G$ is the identity. Then $U^n=U^0$ and the 
scheme is clearly stable. However there are other examples of matrices with 
multiple roots on the unit circle that lead either to bounded or unbounded 
sequences $(G^n)_{n\in\bbN}$. We will call this property $G^n$-boundedness in 
the sequel. It is clearly a property of the amplification matrix and not of 
its characteristic polynomial. If the minimal stable subspaces associated to 
the multiple root are one-dimensional then $G^n$ is bounded (identity example).
If the minimal stable subspaces are multidimensional then $G^n$ grows linearly.
Such cases (which occur for our schemes) should therefore be handled 
specifically.

\subsection{Sketch of proofs}

In the next sections, we will not give the proofs, but only list in a table 
the arguments used for each situation. We describe here the general plan and 
give names to specific final arguments used. The detailed proofs may be found 
in \cite{bidegaray14} for space dimensions 1 and 2. The three dimensional case
is much more tedious and is work in progress.

Usually the system is given in a implicit form. The first step consists in 
writing it in an explicit form. This yields the amplification matrix $G$. Then
we compute its characteristic polynomial $\phi_0$. In order to perform a von 
Neumann analysis, we compute the series $(\phi_m)$. In the general case, under 
the assumption that the stability condition cannot be better than Maxwell's, 
we can apply either Theorem~\ref{Th_Schur} (\textit{Theorem~\ref{Th_Schur}} 
argument) or Theorem~\ref{Th_vonNeumann} (\textit{Theorem~\ref{Th_vonNeumann}} 
argument), check estimates at each level until $\phi_m$ is a one degree 
polynomial. Special cases arise when $\eps_\rms=\eps_\infty$, $\sin(\xi/2)=0$ 
or $\pm1$, and sometimes for limit values of physical coefficients. In these 
cases, different points of view have to be considered:
\begin{itemize}
\item Theorem~\ref{Th_vonNeumann} has to be used instead of 
Theorem~\ref{Th_Schur},
\item Some eigenvalues lie on the unit circle (mostly $\pm1$ or $\pm i$) and 
are simple, it is then sufficient to study only the other eigenvalues 
(\textit{sub-polynomial} argument) and we conclude to a simple von Neumann 
polynomial and stability,
\item Some eigenvalues lie on the unit circle and are not simple, and besides 
the study of the other eigenvalues (to prove that the polynomial is a von 
Neumann one), we have to find out if the associated minimal stable subspaces 
are one- (stable case) or multidimensional (unstable case). This may be 
checked directly on the form of matrix $G$ (\textit{$G$ form} argument), or 
necessitates the computation of eigenvectors (\textit{eigenvectors} argument). 
If only one eigendirection is found for a multiple eigenvalue, the minimal 
subspace is necessarily multidimensional.
\end{itemize}

\section{Debye media}
\label{sec-Debye}

We address two discretizations of Maxwell--Debye equations. The 
first one uses a $(\bfB,\bfE,\bfD)$ setting for the equations and the second a 
$(\bfB,\bfE,\bfP,\bfJ)$ formulation. 


\subsection{Debye--Joseph et al. model}

In \cite{joseph-hagness-taflove}, Joseph et al. close System~\eqref{Max} by a 
discretization for Eq.~\eqref{DebyeD}, namely 
\begin{equation}
\label{DebyeD1}
\begin{array}{l}
\ds \eps_0\eps_\infty t_\rmr \frac{E^{n+1}_j - E^n_j}\dt
+ \eps_0\eps_\rms \frac{E^{n+1}_j + E^n_j}2 
= t_\rmr \frac{D^{n+1}_j - D^n_j}\dt + \frac{D^{n+1}_j + D^n_j}2. 
\end{array}
\end{equation}
System~\eqref{Max}--\eqref{DebyeD1} may be cast in an explicit form which
handles the variable 
\begin{equation*}
U^n_j = \transp(c_\infty B^{n-\frac12}_{j+\frac12}, E^n_j,
D^n_j/\eps_0\eps_\infty)
\end{equation*}
and the amplification matrix $G$ reads  
\begin{equation*}
\left(\begin{array}{ccc}
1 & -\lambda(e^{i\xi}-1) & 0 \\
- \frac{(1+\delta)\lambda(1-e^{-i\xi})}{1+\delta\eps'_\rms} 
& \frac{(1-\delta\eps'_\rms)+(1+\delta)\lambda^2(e^{i\xi}-2+e^{-i\xi})}
{1+\delta\eps'_\rms} 
& \frac{2\delta}{1+\delta\eps'_\rms} \\
-\lambda(1-e^{-i\xi}) & \lambda^2(e^{i\xi}-2+e^{-i\xi}) & 1 
\end{array}\right)
\end{equation*}
where $\lambda=c_\infty\dt/\dx$ is the CFL constant, $\delta=\dt/2t_\rmr>0$ is 
the normalized time step and $\eps'_\rms=\eps_\rms/\eps_\infty\geq1$ denotes 
the normalized static permittivity. Moreover we define 
\begin{equation*}
q=-\lambda^2(e^{i\xi}-2+e^{-i\xi})= 4\lambda^2\sin^2(\xi/2).
\end{equation*}

The characteristic polynomial is proportional to
\begin{eqnarray*}
\phi_0(Z)
& = & [1+\delta\eps'_\rms)]Z^3 - [3 + \delta\eps'_\rms - (1+\delta)q] Z^2 \\
& & + [3-\delta\eps'_\rms - (1-\delta)q] Z - [1-\delta\eps'_\rms].
\end{eqnarray*}
The proofs are summed up in Table~\ref{DJ} and we deduce that the stability 
condition is $q\leq4$ if $\eps_\rms>\eps_\infty$ and $q<4$ if 
$\eps_\rms=\eps_\infty$. 

\begin{table}[ht]
\[
\begin{array}{|c|c|c|c|}
\hline
q & \eps_\rms & \textrm{argument} & \textrm{result} \\
\hline
\hline
]0,4[ & >\eps_\infty 
& \textrm{Theorem~\ref{Th_Schur}} & \textrm{stable} \\
\hline
]0,4[ & =\eps_\infty 
& \textrm{Theorem~\ref{Th_vonNeumann}} & \textrm{stable} \\
\hline
0 & \geq\eps_\infty 
& \textrm{$G$ form} & \textrm{stable} \\
\hline
4 & >\eps_\infty 
& \textrm{Theorem~\ref{Th_vonNeumann}} & \textrm{stable} \\
\hline
4 & =\eps_\infty 
& \textrm{eigenvectors} & \textrm{unstable} \\
\hline
\end{array}
\]
\caption{\label{DJ}Proof arguments and results for the Debye--Joseph et al. 
model.}
\end{table}


\subsection{Debye--Young model}

In \cite{young}, Young closes System~\eqref{MaxP} by two discretizations for 
Eq.~\eqref{DebyeP}, namely 
\begin{equation}
\label{DebyeP1}
t_\rmr \frac{P^{n+\frac12}_j - P^{n-\frac12}_j}\dt 
= - \frac{P^{n+\frac12}_j + P^{n-\frac12}_j}2 
+ \eps_0(\eps_\rms-\eps_\infty) E^n_j, 
\end{equation}
\begin{equation}
\label{DebyeP2}
t_\rmr J^{n+\frac12}_j = - P^{n+\frac12}_j 
+ \eps_0(\eps_\rms-\eps_\infty) \frac{E^{n+1}_j + E^n_j}2.
\end{equation}
Although $J^{n+\frac12}_j$ is used for the computations, this not a genuine
variable for System~\eqref{MaxP}--\eqref{DebyeP1}--\eqref{DebyeP2} which
handles the variable
\begin{equation*}
U^n_j = \transp(c_\infty B^{n-\frac12}_{j+\frac12}, E^n_j,
P^{n-\frac12}_j/\eps_0\eps_\infty)
\end{equation*}
and the amplification matrix $G$ reads  
\begin{equation*}
\left(\begin{array}{ccc}
1 & -\lambda(e^{i\xi}-1) & 0 \\
- \frac{\lambda(1-e^{-i\xi})}{1+\delta\alpha} 
& \frac{1+\delta-\delta\alpha+3\delta^2\alpha -(1+\delta)q}
{(1+\delta)(1+\delta\alpha)} 
& \frac{1-\delta}{1+\delta}\frac{2\delta}{1+\delta\alpha}  \\
0 & \frac{2\delta\alpha}{1+\delta}& \frac{1-\delta}{1+\delta}
\end{array}\right)
\end{equation*}
with the same notation as above and $\alpha=\eps'_\rms-1\geq0$.

The characteristic polynomial is proportional to
\begin{eqnarray*}
\phi_0(Z) 
& = & [(1+\delta\alpha)(1+\delta)] Z^3 
- [3 + \delta + \delta\alpha + 3\delta^2\alpha - (1+\delta)q] Z^2 \\
&& + [3 - \delta - \delta\alpha + 3\delta^2\alpha - (1-\delta)q] Z 
- [(1-\delta\alpha)(1-\delta)].
\end{eqnarray*}
Again, the proofs are summed up in Table~\ref{DY}.

\begin{table}[ht]
\[
\begin{array}{|c|c|c|c|c|}
\hline
q & \eps_\rms & \delta & \textrm{argument} & \textrm{result} \\
\hline
\hline
]0,4] & >\eps_\infty & ]0,1[
& \textrm{Theorem~\ref{Th_Schur}} & \textrm{stable} \\
\hline
]0,4[ & =\eps_\infty & >0
& \textrm{Theorem~\ref{Th_vonNeumann}} & \textrm{stable} \\
\hline
0 & \geq\eps_\infty & >0 
& \textrm{$G$ form} & \textrm{stable} \\
\hline
]0,4] & >\eps_\infty & 1
& \textrm{sub-polynomial} & \textrm{stable} \\
\hline
4 & =\eps_\infty & >0
& \textrm{eigenvectors} & \textrm{unstable} \\
\hline
\end{array}
\]
\caption{\label{DY}Proof arguments and results for the Debye--Young model.}
\end{table}
The stability condition is therefore $q\leq4$ and $\delta\leq1$ if 
$\eps_\rms>\eps_\infty$ and $q<4$ if $\eps_\rms=\eps_\infty$. 

\subsection{Conclusion for one-dimensional Debye schemes}

If $\eps_\rms>\eps_\infty$, the pure CFL condition $q\leq4$ is the same for 
both models. It is exactly the condition for Maxwell equations. However Young 
model necessitates another condition, $\delta\leq1$, which corresponds to a 
sufficient discretization of Debye equation~\eqref{DebyeP}. Even if we are 
interested here in stability properties, such conditions are to be taken to 
ensure equations to be correctly taken into account. Results are given in 
physical variables in Table~\ref{Debyeresume}.

\begin{table}[ht]
\[
\begin{array}{|c|c|c|}
\hline
\textrm{Scheme} & & \textrm{dimension 1} \\
\hline\hline
\multicolumn{3}{|c|}{\eps_\rms>\eps_\infty} \\
\hline\hline
\textrm{Joseph et al.} & q\leq4 & \dt\leq\frac{\dx}{c_\infty} \\
\hline
\textrm{Young} & q\leq4,\ \delta\leq1 
& \dt\leq\min(\frac{\dx}{c_\infty},2t_\rmr) \\
\hline\hline
\multicolumn{3}{|c|}{\eps_\rms=\eps_\infty} \\
\hline\hline
\textrm{Joseph et al.} & q<4 & \dt<\frac{\dx}{c_\infty} \\
\hline
\textrm{Young} & q<4 & \dt<\frac{\dx}{c_\infty} \\
\hline
\end{array}
\]
\caption{\label{Debyeresume} Stability of Debye models for 
$\eps_\rms>\eps_\infty$ and $\eps_\rms=\eps_\infty$.}
\end{table}

To compare conditions on $q$ and $\delta$, let us consider a simple physical 
case. We assume that a matter with $\eps_\infty=1$ (and thus
$c_\infty\simeq 3\,10^{8}\uv$) is lighted by an optical wave of say 
wavelength 1\,$\mu$m. The space step $\dx$ has to be smaller than this 
wavelength, and therefore $q<4$ reads at least $\dt<\frac13\,10^{-14}\ut$. In 
a Debye medium, relaxation times $t_\rmr$ are of the order of a picosecond (or 
even a nanosecond) which is many decades larger than the previous bound. The 
estimate $q<4$ is thus predominant and both models present the same advantages.
Only the value of $\eps_\infty$ yields the CFL condition. A typical example is
water for which $\eps_\infty=1.8$, $\eps_\rms=81.0$ and 
$t_\rmr=9.4\,10^{-12}\ut$ \cite{young-kittichartphayak-kwok-sullivan}. 
Condition $\dt\leq2t_\rmr$ comes to $\dt\leq 1.88\,10^{-11}\ut$. Condition 
$q\leq4$ yields a similar condition if $\dx=4.2\,10^{-3}\um$. This is of 
course much larger than any reasonable space step for Maxwell equations and 
optical waves. The stability condition for water is $q<4$ for both schemes. A 
quite different material is for example the 0.25-dB loaded foam given in 
\cite{luebbers-steich-kunz} for which $\eps_\infty=1.01$, $\eps_\rms=1.16$ and 
$t_\rmr=6.497\,10^{-10}\ut$. Condition $\dt\leq2t_\rmr$ comes to 
$\dt\leq 1.3\,10^{-9}\ut$ and $q\leq4$ yields a similar condition if 
$\dx=3.9\,10^{-1}\um$. Once more, the stability condition for water is $q<4$ 
for both schemes.

In conclusion for current material the stability condition is the same for 
Maxwell--Debye equations as for the usual Yee scheme. The result announced in
\cite{petropoulos1} was $q\leq4$ for Joseph et al. scheme and for water, which
is consistent with our result.

\section{Lorentz media}
\label{sec-Lorentz}

Three discretizations of Maxwell--Lorentz equations are now 
addressed. The first one uses a $(\bfB,\bfE,\bfD)$ setting and the two others 
a $(\bfB,\bfE,\bfP,\bfJ)$ formulation, but differ from the time-discretization 
of $\bfJ$. 

Each of these models reads the same in the harmonic ($\nu=0$) or an-harmonic 
($\nu>0$) cases. However the analysis will differ greatly since $\phi_1\equiv0$
for all the schemes in the harmonic cases. 


\subsection{Lorentz--Joseph et al. model}

In \cite{joseph-hagness-taflove}, system~\eqref{Max} is closed by a 
discretization for Eq.~\eqref{DebyeD}, namely 
\begin{equation}
\label{LorentzD1}
\begin{array}{l}
\ds \eps_0\eps_\infty \frac{E^{n+1}_j - 2E^n_j + E^{n-1}_j}{\dt^2}
+ \nu \eps_0\eps_\infty \frac{E^{n+1}_j - E^{n-1}_j}{2\dt} 
+ \eps_0\eps_\rms \omega_1^2 \frac{E^{n+1}_j + E^{n-1}_j}2 \\ 
\ds \hspace{5mm} 
= \frac{D^{n+1}_j - 2D^n_j + D^{n-1}_j}{\dt^2}
+ \nu \frac{D^{n+1}_j - D^{n-1}_j}{2\dt}  
+ \omega_1^2 \frac{D^{n+1}_j + D^{n-1}_j}2 \\
\end{array}
\end{equation}
The explicit version of system~\eqref{Max}--\eqref{LorentzD1} does not use 
explicitly the value of $D^{n-1}_j$ and therefore this system handles the 
variable 
\begin{equation*}
U^n_j = \transp(c_\infty B^{n-\frac12}_{j+\frac12}, E^n_j, E^{n-1}_j,
D^n_j/\eps_0\eps_\infty).
\end{equation*}
The amplification matrix $G$ reads  
\begin{equation*}
\left(\begin{array}{cccc}
1 & -\lambda(e^{i\xi}-1) & 0 & 0 \\
- \frac{2\delta\lambda(1-e^{-i\xi})}{1+\delta+\omega\eps'_\rms} 
& \frac{2-q(1+\delta+\omega)}{1+\delta+\omega\eps'_\rms} 
& \frac{1-\delta+\omega\eps'_\rms}{1+\delta+\omega\eps'_\rms} 
& \frac{2\omega}{1+\delta+\omega\eps'_\rms} \\
0 & 1 & 0 & 0 \\
- \lambda(1-e^{-i\xi}) & -q & 0 & 1
\end{array}\right)
\end{equation*}
where $\delta=\nu\dt/2\geq0$ is the new normalized time step, and 
$\omega=\omega_1^2\dt^2/2>0$ denotes the normalized squared frequency. The 
other notations used for the Debye model remain valid. 

The characteristic polynomial is proportional to
\begin{eqnarray*}
\phi_0(Z) 
& = & [1+\delta+\omega\eps'_\rms] Z^4 
- [4 + 2\delta + 2\omega\eps'_\rms - (1 + \delta + \omega) q] Z^3 \\
&&+ [6 + 2\omega\eps'_\rms - 2q] Z^2 
- [4 - 2\delta + 2\omega\eps'_\rms - (1-\delta+\omega) q] Z \\
&& + [1-\delta+\omega\eps'_\rms].
\end{eqnarray*}
The proofs are summed up in Table~\ref{LJ} for the an-harmonic and the 
harmonic case.
\begin{table}[ht]
\[
\begin{array}{|c|c|c|c|}
\hline
q & \eps_\rms & \textrm{argument} & \textrm{result} \\
\hline
\hline
\multicolumn{4}{|c|}{\textrm{an-harmonic: }\nu>0} \\
\hline
\hline
]0,2[ & >\eps_\infty 
& \textrm{Theorem~\ref{Th_Schur}} & \textrm{stable} \\
\hline
]0,2] & =\eps_\infty 
& \textrm{Theorem~\ref{Th_vonNeumann}} & \textrm{stable} \\
\hline
0 & \geq\eps_\infty 
& \textrm{$G$ form} & \textrm{stable} \\
\hline
2 & \geq\eps_\infty 
& \textrm{sub-polynomial} & \textrm{stable} \\
\hline
\hline
\multicolumn{4}{|c|}{\textrm{harmonic: }\nu=0} \\
\hline
\hline
]0,2[ & >\eps_\infty 
& \textrm{Theorem~\ref{Th_vonNeumann}} & \textrm{stable} \\
\hline
]0,2] & =\eps_\infty 
& \textrm{sub-polynomial} & \textrm{unstable} \\
\hline
0 & \geq\eps_\infty 
& \textrm{$G$ form} & \textrm{stable} \\
\hline
2 & \geq\eps_\infty 
& \textrm{sub-polynomial} & \textrm{stable} \\
\hline
\end{array}
\]
\caption{\label{LJ}Proof arguments and results for the Lorentz--Joseph et al. 
model.}
\end{table}

In the an-harmonic case the stability condition is $q\leq2$ whatever 
$\eps_\rms\geq\eps_\infty$ is. The $\eps_\rms=\eps_\infty$ harmonic case, 
needs some explanation. For $q\in]0,2]$, $\phi_0$ may be cast as the product 
of two second order polynomials. The roots are two couples of conjugate 
complex roots of modulus 1. For the specific value $q=2\omega/(1+\omega)$, 
which always lies in the interval $]0,2]$, the two couples degenerate in one
double couple, and the associated minimal stable sub-spaces are 
two-dimensional. To avoid this instability one may think to bound $q$ and say
that the scheme is stable provided $q\in[0,2\omega/(1+\omega)[$. But if we 
come back to the original variables, we see that this is not an upper bound on 
$\dt$ but rather a lower bound on $\dx$, which we surely do not want. It is 
therefore better to avoid using Joseph et al. scheme in this very specific 
case, $\eps_\rms=\eps_\infty$ and $\nu=0$, and we hope to find a better scheme 
for this case in the following examples.


\subsection{Lorentz--Kashiwa et al. model}

In \cite{kashiwa-yoshida-fukai}, Kashiwa et al. close a modified version of 
System~\eqref{MaxP}, which consists of the three first equations in 
System~\eqref{LorentzP1}, by a discretization for Eq.~\eqref{LorentzP}, namely 
\begin{equation}
\label{LorentzP1}
\begin{array}{rcl}
\ds \frac1\dt (B_{j+\frac12}^{n+\frac12}-B_{j+\frac12}^{n-\frac12})
& = & \ds - \frac1\dx (E_{j+1}^n-E_j^n), \\
\ds \frac{\eps_0\eps_\infty}\dt (E_j^{n+1}-E_j^n)  
& = & \ds 
- \frac1{\mu_0 \dx} (B_{j+\frac12}^{n+\frac12}-B_{j-\frac12}^{n+\frac12}) 
- \frac1\dt (P^{n+1}_j-P^n_j), \\
\ds \frac1\dt (P^{n+1}_j-P^n_j) & = & \ds \frac12 (J^{n+1}_j+J^n_j), \\
\ds \frac1\dt (J^{n+1}_j-J^n_j) & = & \ds - \frac\nu2 (J^{n+1}_j+J^n_j) 
+ \frac{\omega_1^2(\eps_\rms-\eps_\infty)\eps_0}2 (E^{n+1}_j+E^n_j) \\
&& - \frac{\omega_1^2}2 (P^{n+1}_j+P^n_j).
\end{array}
\end{equation}
The explicit version of system~\eqref{LorentzP1} handles the variable 
\begin{equation*}
U^n_j = \transp(c_\infty B^{n-\frac12}_{j+\frac12}, E^n_j, 
P^n_j/\eps_0\eps_\infty,\dt J^n_j/\eps_0\eps_\infty)
\end{equation*}
and the amplification matrix $G$ reads  
\begin{equation*}
\left(\begin{array}{cccc}
1 & -\lambda(e^{i\xi}-1) & 0 & 0 \\
\frac{-\lambda(1-e^{-i\xi})(\Delta-\frac12\omega\alpha)}{\Delta} 
& \frac{\Delta-q\Delta-(2-q)\frac12\omega\alpha}{\Delta} 
& \frac{\omega}{\Delta} & \frac{-1}{\Delta} \\
\frac{-\lambda(1-e^{-i\xi})\frac12\omega\alpha}{\Delta} 
& \frac{(2-q)\frac12\omega\alpha}{\Delta} 
& \frac{\Delta-\omega}{\Delta} & \frac{1}{\Delta} \\
\frac{-\lambda(1-e^{-i\xi})\omega\alpha}{\Delta} 
& \frac{(2-q)\omega\alpha}{\Delta} & \frac{-2\omega}{\Delta} 
& \frac{2-\Delta}{\Delta}
\end{array}\right)
\end{equation*}
where together with the previously defined notations, 
$\Delta=1+\delta+\omega\eps'_\rms/2$. 

The characteristic polynomial is proportional to
\begin{eqnarray*}
\phi_0(Z)
& = & [1+\delta+\frac12\omega\eps'_\rms] Z^4 
- [4+2\delta-(1+\delta+\frac12\omega)q] Z^3 \\
&& + [6-\omega\eps'_\rms+(\omega-2)q] Z^2 
- [4-2\delta-(1-\delta+\frac12\omega)q] Z \\
&& + [1-\delta+\frac12\omega\eps'_\rms].
\end{eqnarray*}
The proofs are summed up in Table~\ref{LK}. Both in the an-harmonic and 
harmonic cases, the stability condition is $q<4$ which is much better than the 
previous scheme since we gain a factor 2 on $\dt$ and we have no problem when 
$\eps_\rms=\eps_\infty$ and $\nu=0$ as for the previous model.

\begin{table}[ht]
\[
\begin{array}{|c|c|c|c|}
\hline
q & \eps_\rms & \textrm{argument} & \textrm{result} \\
\hline
\hline
\multicolumn{4}{|c|}{\textrm{an-harmonic: }\nu>0} \\
\hline
\hline
]0,4[ & >\eps_\infty & \textrm{Theorem~\ref{Th_Schur}} & \textrm{stable}\\
\hline
]0,4[ & =\eps_\infty & \textrm{Theorem~\ref{Th_vonNeumann}} 
& \textrm{stable}\\
\hline
0 & \geq\eps_\infty & \textrm{G form} & \textrm{stable}\\
\hline
4 & \geq\eps_\infty & \textrm{eigenvectors} & \textrm{unstable}\\
\hline
\hline
\multicolumn{4}{|c|}{\textrm{harmonic: }\nu=0} \\
\hline
\hline
]0,4[ & \geq\eps_\infty & \textrm{Theorem~\ref{Th_vonNeumann}} 
& \textrm{stable}\\
\hline
0 & \geq\eps_\infty & \textrm{G form} & \textrm{stable}\\
\hline
4 & \geq\eps_\infty & \textrm{eigenvectors} & \textrm{unstable}\\
\hline
\end{array}
\]
\caption{\label{LK}Proof arguments and results for the Lorentz--Kashiwa et 
al. model.}
\end{table}


\subsection{Lorentz--Young model}

In \cite{young}, System~\eqref{MaxP} is closed by a discretization for 
Eq.~\eqref{LorentzP}, namely 
\begin{equation}
\label{LorentzP2}
\begin{array}{l}
\ds \frac1\dt (P^{n+1}_j-P^n_j) = J^{n+\frac12}, \\
\ds \frac1\dt (J^{n+\frac12}_j-J^{n-\frac12}_j)
= - \frac\nu2 (J^{n+\frac12}_j+J^{n-\frac12}_j) \\
\ds \hspace{3cm}
+ \omega_1^2(\eps_\rms-\eps_\infty)\eps_0 E^n_j - \omega_1^2 P^n_j.
\end{array}
\end{equation}
The explicit version of System~\eqref{MaxP}--\eqref{LorentzP2} handles once 
more the variable 
\begin{equation*}
U^n_j = \transp(c_\infty B^{n-\frac12}_{j+\frac12}, E^n_j, 
P^n_j/\eps_0\eps_\infty,\dt J^n_j/\eps_0\eps_\infty)
\end{equation*}
and the amplification matrix $G$ reads  
\begin{equation*}
\left(\begin{array}{cccc}
1 & -\lambda(e^{i\xi}-1) & 0 & 0 \\
-\lambda(1-e^{-i\xi}) & \frac{(1-q)(1+\delta)-2\omega\alpha}{1+\delta} 
& \frac{2\omega}{1+\delta} & -\frac{1-\delta}{1+\delta} \\
0 & \frac{2\omega\alpha}{1+\delta} 
& \frac{1+\delta-2\omega}{1+\delta} & \frac{1-\delta}{1+\delta} \\
0 & \frac{2\omega\alpha}{1+\delta} 
& \frac{-2\omega}{1+\delta} & \frac{1-\delta}{1+\delta} \\
\end{array}\right)
\end{equation*}

The characteristic polynomial is proportional to
\begin{eqnarray*}
\phi_0(Z)
& = & [1+\delta] Z^4 - [4 + 2\delta - 2\omega\eps'_\rms - (1+\delta)q] Z^3 \\
&& + 2 [3 - 2\omega\eps'_\rms + (\omega-1)q ] Z^2 \\
&& - [4 - 2\delta - 2\omega\eps'_\rms - (1-\delta)q] Z + [1-\delta].
\end{eqnarray*}
The proofs are summed up in Table~\ref{LY}. This scheme combines three 
drawbacks we have already encountered. First as for the Debye model, there 
is an extra condition on the time step: $\omega<2/(2\eps'_\rms-1)$. This will
have to be compared to the condition on $q$ for physical examples. Second, as 
for the Lorentz--Joseph et al. scheme we need a twice smaller $\dt$ than for
raw Maxwell equations: $q\leq2$ instead of $q\leq4$. Last, and also as for the 
Lorentz-Joseph et al. model, the $\eps_\rms=\eps_\infty$ and $\nu=0$ leads to 
an instability. This is exactly the same story. This time $q=2\omega$ leads to
double couples of conjugate complex roots of modulus 1, with two-dimensional
minimal stable sub-spaces. If $\omega<1$ this value of $q$ is however never 
reached, but $\omega<1$ is a stronger assumption than 
$\omega<2/(2\eps'_\rms-1)$. We will see what this amounts to in numerical 
applications.
 
\begin{table}[ht]
\[
\begin{array}{|c|c|c|c|c|}
\hline
q & \eps_\rms & \omega & \textrm{argument} & \textrm{result} \\
\hline
\hline
\multicolumn{5}{|c|}{\textrm{an-harmonic: }\nu>0} \\
\hline
\hline
]0,2[ & > \eps_\infty & \leq \frac{2}{2\eps'_\rms-1} 
& \textrm{Theorem~\ref{Th_Schur}} & \textrm{stable} \\
\cline{1-3}
2 & > \eps_\infty & < \frac{2}{2\eps'_\rms-1} & & \\
\hline
]0,2] & = \eps_\infty & < 2 & \textrm{Theorem~\ref{Th_vonNeumann}} 
& \textrm{stable} \\
\hline
]0,2] & = \eps_\infty & = 2 & \textrm{sub-polynomial} & \textrm{stable} \\
\hline
2 & > \eps_\infty & = \frac{2}{2\eps'_\rms-1}
 & \textrm{Theorem~\ref{Th_vonNeumann}} & \textrm{stable} \\
\hline
0 & \geq \eps_\infty & \leq \frac{2}{2\eps'_\rms-1} & \textrm{G form} 
& \textrm{stable} \\
\hline
\hline
\multicolumn{5}{|c|}{\textrm{harmonic: }\nu=0} \\
\hline
\hline
]0,2[ & > \eps_\infty & \leq \frac{2}{2\eps'_\rms-1} 
& \textrm{Theorem~\ref{Th_vonNeumann}} & \textrm{stable} \\
\cline{1-3}
2 & > \eps_\infty & < \frac{2}{2\eps'_\rms-1} & & \\
\hline
]0,2] & = \eps_\infty & < 2 & \textrm{eigenvectors} 
& \textrm{unstable} \\
\hline
]0,2] & = \eps_\infty & = 2 & \textrm{Theorem~\ref{Th_vonNeumann}} 
& \textrm{stable} \\
\hline
2 & > \eps_\infty & = \frac{2}{2\eps'_\rms-1}
 & \textrm{eigenvectors} & \textrm{unstable} \\
\hline
0 & > \eps_\infty & \leq \frac{2}{2\eps'_\rms-1} & \textrm{G form} 
& \textrm{stable} \\
\cline{1-3}
0 & = \eps_\infty & < \frac{2}{2\eps'_\rms-1} & & \\
\hline
0 & = \eps_\infty & = \frac{2}{2\eps'_\rms-1} & \textrm{eigenvectors} 
& \textrm{unstable} \\
\hline
\end{array}
\]
\caption{\label{LY}Proof arguments and results for the Lorentz--Young model.}
\end{table}

\subsection{Conclusion for one-dimensional Lorentz schemes}
 
We can summarize all our results for Lorentz schemes in 
Table~\ref{Lorentzresume}. We chose not to translate the result for the Young 
scheme for $\eps_\rms=\eps_\infty$ as a condition on $\dx$ ($q<2\omega$) but 
as a condition on $\dt$ ($\omega<1$, and therefore $q=2\omega$ is not 
reached).\\
For the harmonic Young scheme if $\eps_\rms>\eps_\infty$ the condition is 
slightly better since $q=2$ and $\omega<2/(2\eps'_\rms-1)$, or $q<2$ and 
$\omega=2/(2\eps'_\rms-1)$ also yield stable schemes. \\

Contrarily to Debye materials, for which Joseph et al. model and Young model 
compete, the Kashiwa et al. model seems to overcome others for Lorentz 
material. First, there is a gain in CFL condition $q<4$ is twice better as 
$q\leq2$, second, there are no instabilities for limiting values of the 
physical coefficients and last there are no extra condition on the time step. 
In practice, an extra condition is however needed to account for the dynamics 
of the Lorentz equation, but not for stability reasons.

However we can compare the relative strength of the different conditions on 
$\dt$ for Joseph et al. and Young models. The values used in 
\cite{petropoulos1} are $\eps_\infty=1$, $\eps_\rms=2.25$, 
$\omega_1=4\,10^{16}\uf$ and $\nu=0.56\,10^{16}\uf$. Condition 
$\omega\leq2/\sqrt{2\eps'_\rms-1}$ comes to $\dt \leq 2.7\,10^{-17}\ut$ which 
is very small and corresponds to $\dx=1.13\,1^{-8}\um$ in the $q<2$ condition. 
This space step is more than sufficient to discretize optical waves. For such 
a material the extra condition imposed by the Joseph et al. scheme is stronger
than the basic CFL condition. The Kashiwa et al. model is then more advisable. 

\begin{table}[ht]
\[
\begin{array}{|c|c|c|}
\hline
\textrm{Scheme} & & \textrm{dimension 1} \\
\hline\hline
\multicolumn{3}{|c|}{\textrm{an-harmonic: }\nu>0, 
\textrm{ and } \eps_\rms\geq\eps_\infty} \\
\hline\hline
\textrm{Joseph} & q\leq2 
& \dt\leq\frac{\dx}{\sqrt2 c_\infty} \\
\hline
\textrm{Kashiwa} & q<4 
& \dt<\frac{\dx}{c_\infty} \\
\hline
\textrm{Young} 
& \begin{array}{c}q\leq2, \\ \omega\leq\frac2{2\eps'_\rms-1} \end{array}
& \dt\leq\min(\frac{\dx}{\sqrt2 c_\infty},
\frac2{\omega_1\sqrt{2\eps'_\rms-1}})\\ 
\hline\hline
\multicolumn{3}{|c|}{\textrm{harmonic: }\nu=0, 
\textrm{ and } \eps_\rms>\eps_\infty} \\
\hline\hline
\textrm{Joseph} & q\leq2 
& \dt\leq\frac{\dx}{\sqrt2 c_\infty} \\
\hline
\textrm{Kashiwa} & q<4 
& \dt<\frac{\dx}{c_\infty} \\
\hline
\textrm{Young} 
& \begin{array}{c}q<2, \\ \omega<\frac2{2\eps'_\rms-1} \end{array}
& \dt<\min(\frac{\dx}{\sqrt2 c_\infty},\frac2{\omega_1\sqrt{2\eps'_\rms-1}})\\ 
\hline\hline
\multicolumn{3}{|c|}{\textrm{harmonic: }\nu=0, 
\textrm{ and } \eps_\rms=\eps_\infty} \\
\hline\hline
\textrm{Joseph} & q<\frac{2\omega}{1+\omega} 
& \textrm{condition on }\dx \\
\hline
\textrm{Kashiwa} & q<4 
& \dt<\frac{\dx}{c_\infty} \\
\hline
\textrm{Young} 
& \begin{array}{c}q<2,\\\omega<1 \end{array}
& \dt<\min(\frac{\dx}{\sqrt2 c_\infty},\frac{\sqrt2}{\omega_1})\\ 
\hline
\end{array}
\]
\caption{\label{Lorentzresume}Stability of an-harmonic and harmonic Lorentz 
models for $\eps_\rms>\eps_\infty$ and $\eps_\rms>\eps_\infty$.}
\end{table}

In \cite{young} there is a totally different material for which 
$\eps_\infty=1.5$, $\eps_\rms=3$, $\omega_1=2\pi\,5\,10^{10}\uf$ and 
$\nu=10^{10}\uf$ (these round values certainly refer to a model material). In 
this case $\omega\leq2/\sqrt{2\eps'_\rms-1}$ comes to $\dt\leq3.6\,10^{-12}\ut$
which corresponds to $\dx=1.9\,1^{-3}\um$ in the $q<2$ condition. For this 
material condition $q<2$ is the strongest for optical waves. The Kashiwa et al.
model is however more advisable, since it allows $q<4$ instead of $q\leq2$. 

The results obtained in \cite{petropoulos1} where obtained for our first cited 
material and for Joseph et al. and Kashiwa et al. models. He observed 
instabilities for $\xi>\frac\pi2$. We note that if $\xi\leq\frac\pi2$ then 
$\sin(\xi/2)\leq1/\sqrt{2}$ and $q\leq2$ instead of $q\leq4$. This is exactly 
our result. He found also the Kashiwa et al. scheme to stable for $q\leq4$.
 
\section{Two-dimensional results}
\label{sec-bidim}

In a two-dimensional context where unknowns depend only on space
variables $x$ and $y$, Maxwell system may be split in two decoupled systems 
corresponding to the transverse electric (TE) ($B_x$, $B_y$, $E_z$) and the 
transverse magnetic (TM) ($B_z$, $E_x$, $E_y$) polarizations. In the 
one-dimensional case, Maxwell--Debye equations were represented by three 
equations and Maxwell--Lorentz by four equations. In the TE polarization, one
more Faraday equation is added and we have four equations for Maxwell--Debye 
and five equations for  Maxwell--Lorentz. In the TM polarization for the 
Maxwell--Debye model, one Amp\`ere equation and one Debye equation have to be 
added, leading to five equations systems. For the Maxwell--Lorentz model, 
there are one Amp\`ere equation and two Lorentz equations more, and the system
consists of seven equations. 

The principle of the stability analysis is exactly the same, but we now have
larger polynomials to study. A small miracle however happens: one-dimensional
polynomials are a factor in two-dimensional polynomials. More precisely we 
now denote by $\dx_x$ and $\dx_y$ the space steps in the $x$- and 
$y$-directions respectively and by $q$ the quantity 
\begin{equation*}
q = q_x + q_y 
= 4c_\infty^2 \left(\frac{\dt^2}{\dx_x^2}\sin^2(\xi_x/2)
               +\frac{\dt^2}{\dx_y^2}\sin^2(\xi_y/2)\right)
\end{equation*}
(recall $q=4c_\infty^2\frac{\dt^2}{\dx_x^2}\sin^2(\xi_x/2)$ in 1D). Then in 
the two-dimensional TE polarization 
\begin{equation*}
\phi^{2D,TE}_0(Z) = [Z-1]\phi^{1D}_0(Z),
\end{equation*}
for all the Maxwell--Debye and Maxwell-Lorentz schemes we study here. This 
could be a problem, if 1 is already a root of $\phi^{1D}_0(Z)$, i.e. when 
$q=0$, but it happens that it is never a problem: minimal stable sub-spaces are
always one-dimensional. In the TM polarization, the same factorization occurs
but the remaining polynomial is slightly more complicated, namely 
\begin{equation*}
\phi^{2D,TM}_0(Z) = [Z-1]\psi_0(Z)\phi^{1D}_0(Z),
\end{equation*}
where $\psi_0(Z)$ 
is equal to: \\
-- Debye--Joseph et al. model
\begin{equation*}
[(1+\delta\eps'_\rms)Z-(1-\delta\eps'_\rms)].
\end{equation*} 
-- Debye--Young model
\begin{equation*}
[(1+\alpha)(1+\delta\alpha)Z-(1-\alpha)(1-\delta\alpha)].
\end{equation*} 
-- Lorentz--Joseph et al. model
\begin{equation*}
[(1+\delta+\omega\eps'_\rms)Z^2-2Z+(1-\delta+\omega\eps'_\rms)].
\end{equation*} 
- Lorentz--Kashiwa et al. model
\begin{equation*}
[(1+\delta+\frac12\omega\eps'_\rms)Z^2-(2-\omega\eps'_\rms)Z
+(1-\delta+\frac12\omega\eps'_\rms)].
\end{equation*} 
-- Lorentz--Young model
\begin{equation*}
[(1+\delta)Z^2-2(1-\omega\eps'_\rms)Z+(1-\delta)].
\end{equation*} 
As for the TE polarization the extra eigenvalue 1 is never a source of 
instability. The other extra eigenvalues always lie inside or on the unit 
circle (conjugate complex roots). The only problem is when modulus 1 
eigenvalues are also eigenvalues of the one-dimensional polynomial. This only
occurs for the Lorentz-Joseph et al. scheme is $\eps_\rms=\eps_\infty$, and 
$q=2\omega/(1+\omega)$, which is a resonant value we have already pointed out 
in the harmonic case for this scheme. 

We shall not duplicate Tables~\ref{Debyeresume} and \ref{Lorentzresume} for 
two-dimensional models. If $\dx_x=\dx_y\equiv\dx$, condition $q\leq4$ becomes 
$\dt\leq\dx/(\sqrt{2}c_\infty)$ and condition $q\leq2$ becomes 
$\dt\leq\dx/(2c_\infty)$ in the physical variables. Besides, Lorentz--Joseph 
et al. model which was leading to a lower bound on $\dx$ in the harmonic case, 
leads also to such a bound in the an-harmonic case. These are the only 
differences with Tables~\ref{Debyeresume} and \ref{Lorentzresume}.

\section{Conclusion}

We have studied a class of FD--TD schemes for dispersive materials based on 
the Yee scheme for Maxwell equations and compared them from the stability 
point of view. This study was inspired by Petropoulos \cite{petropoulos1} 
who performs the same analysis but using specific values for the physical and 
numerical constants and using numeric routines to locate eigenvalues of the 
amplification matrix. Here we have general results which gives you the 
constraint on numerical constants ($\dt$ and $\dx$) for any Debye or Lorentz 
material. Our results confirm those of Petropoulos.

For usual Debye media, both studied schemes are stable under the same
conditions as the Yee scheme, ensuring also, if applied to optical waves, a
fine discretization of the Debye equation. Among the studied schemes for 
Lorentz media, Kashiwa et al. model clearly ranks first as far as stability is 
concerned., Its stability condition is also that of the Yee scheme. However to 
take properly into account the Lorentz model, a smaller time step may have to 
be chosen, independently of stability issues. Such results have been proved 
for 1D and 2D models. The 3D case, which is much more tedious, is being 
studied and analogous results are expected.


\begin{thebibliography}{10}

\bibitem{petropoulos1}
P.G. Petropoulos,
\newblock ``Stability and phase error analysis of {FD--TD} in dispersive
  dielectrics,''
\newblock {\em {IEEE} {T}ransactions on {A}ntennas and {P}ropagation}, vol. 42,
  no. 1, pp. 62--69, 1994.

\bibitem{luebbers-hunsberger-kunz-standler-schneider}
R.~Luebbers, F.P. Hunsberger, K.S. Kunz, R.B. Standler, and M.~Schneider,
\newblock ``A frequency-dependent finite-difference time--domain formulation
  for dispersive materials,''
\newblock {\em {IEEE} {T}ransactions on {E}lectromagnetic {C}ompatibility},
  vol. 32, no. 3, pp. 222--227, 1990.

\bibitem{young-kittichartphayak-kwok-sullivan}
J.L. Young, A.~Kittichartphayak, Y.M. Kwok, and D.~Sullivan,
\newblock ``On the dispersion errors related to {FD$^2$TD} type schemes,''
\newblock {\em {IEEE} {T}ransactions on {A}ntennas and {P}ropagation}, vol. 43,
  no. 8, pp. 1902--1910, 1995.

\bibitem{yee}
K.S. Yee,
\newblock ``Numerical solution of initial boundary value problems involving
  {M}axwell's equations in isotropic media,''
\newblock {\em {IEEE} {T}ransactions on {A}ntennas and {P}ropagation}, vol. 14,
  no. 3, pp. 302--307, 1966.

\bibitem{feise-schneider-bevelacqua}
M.W. Feise, J.B. Schneider, and P.J. Bevelacqua,
\newblock ``Finite-difference and pseudospactral time--domain methods applied
  to backward-wave metamaterials,''
\newblock {\em {IEEE} {T}ransactions on {A}ntennas and {P}ropagation}, vol. 52,
  no. 11, pp. 2955--2962, 2004.

\bibitem{stoykov-kuiken-lowery-taflove}
N.S. Stoykov, T.A. Kuiken, M.M. Lowery, and A.~Taflove,
\newblock ``Finite-element time--doamin algorithms for modeling linear {D}ebye
  and {L}orentz dielectric dispersions at low frequencies,''
\newblock {\em {IEEE} {T}ransactions on {B}iomedical {E}ngineering}, vol. 50,
  no. 9, pp. 1100--1107, 2003.

\bibitem{kashiwa-yoshida-fukai}
T.~Kashiwa, N.~Yoshida, and I.~Fukai,
\newblock ``A treatment by the {FD--TD} method of the dispersive
  characteristics associated with orientation polarization,''
\newblock {\em {I}nstitute of {E}lectronics, {I}nformation and {C}ommunication
  {E}ngineers {T}ransactions}, vol. E73, pp. 1326--1328, 1990.

\bibitem{joseph-hagness-taflove}
R.M. Joseph, S.C. Hagness, and A.~Taflove,
\newblock ``Direct time integration of {M}axwell's equations in linear
  dispersive media with absorption for scattering and propagation of
  femtosecond electromagnetic pulses,''
\newblock {\em {O}ptical {L}etters}, vol. 16, no. 18, pp. 1412--1414, 1991.

\bibitem{young}
J.L. Young,
\newblock ``Propagation in linear dispersive media: {F}inite difference
  time--domain methodologies,''
\newblock {\em {IEEE} {T}ransactions on {A}ntennas and {P}ropagation}, vol. 43,
  no. 4, pp. 422--426, 1995.

\bibitem{strikwerda}
J.C. Strikwerda,
\newblock {\em Finite Difference Schemes and Partial Differential Equations},
\newblock Wadworth \& Brooks/Cole, 1989.

\bibitem{bidegaray14}
B.~Bid\'egaray-Fesquet,
\newblock ``Analyse de von {N}eumann de sch\'emas aux diff\'erences finies pour
  les \'equations de {M}axwell--{D}ebye et de {M}axwell--{L}orentz,''
\newblock Tech. {R}ep., LMC-IMAG, 2005.


\bibitem{luebbers-steich-kunz}
R.~Luebbers, D.~Steich, and K.~Kunz,
\newblock ``{FDTD} calculation of scattering from frequency-dependent
  materials,''
\newblock {\em {IEEE} {T}ransactions on {A}ntennas and {P}ropagation}, vol. 41,
  no. 9, pp. 1249--1257, 1993.

\end{thebibliography}
\end{document}